\documentclass{article}
\usepackage{amssymb}
\usepackage{amsfonts}
\usepackage{amsmath}

\input{tcilatex}
\begin{document}

\title{On the time evolution of Bernstein processes associated with a class
of parabolic equations }
\author{Pierre-A. Vuillermot \\
Inst. \'{E}lie Cartan de Lorraine, UMR-CNRS 7502, Nancy, France}
\date{}
\maketitle

\begin{abstract}
In this article dedicated to the memory of Igor D. Chueshov, I first
summarize in a few words the joint results that we obtained over a period of
six years regarding the long-time behavior of solutions to a class of
semilinear stochastic parabolic partial differential equations. Then, as the
beautiful interplay between partial differential equations and probability
theory always was close to Igor's heart, I present some new results
concerning the time evolution of certain Markovian Bernstein processes
naturally associated with a class of deterministic linear parabolic partial
differential equations. Particular instances of such processes are certain
conditioned Ornstein-Uhlenbeck processes, generalizations of Bernstein
bridges and Bernstein loops, whose laws may evolve in space in a non trivial
way. Specifically, I examine in detail the time development of the
probability of finding such processes within two-dimensional geometric
shapes exhibiting spherical symmetry. I also define a Faedo-Galerkin scheme
whose ultimate goal is to allow approximate computations with controlled
error terms of the various probability distributions involved.
\end{abstract}

\section{Introduction and outline}

This article is a tribute to some of the works and achievements of our
friend and colleague Igor D. Chueshov, who unfortunately and unexpectedly
passed away on April 23rd, 2016. The qualitative analysis of the behavior of
solutions to various stochastic partial differential equations, henceforth
SPDEs, was one of Igor's strong points. I have therefore deemed it
appropriate to briefly summarize here the results that he and I obtained in
that area over a period stretching from 1998 to 2004. As far as the
presentation of the many other facets of his activities is concerned, I am
thus referring the reader to the other contributions in this volume.

When Igor and I first met in 1994 on the occasion of an international
conference on SPDEs in Luminy, we set out to investigate the behavior of
solutions to those stochastic parabolic equations which specifically occur
in population dynamics, population genetics, nerve pulse propagation and
related topics, given the fact that there were already a substantial number
of works in those areas concerning the deterministic case (see, e.g., \cite%
{bernfeldhuvui}, \cite{vuillermot} and the many references therein). But
instead of starting up front with partial differential equations driven by
some kind of noise, we first considered a class of random parabolic
initial-boundary value problems mainly for the sake of simplification.
Assuming then various statistical and dynamical properties such as those of
the central limit theorem and the Ornstein-Uhlenbeck process for the
lower-order coefficients of the equations, we eventually elucidated the
ultimate behavior of the corresponding solution random fields in \cite%
{chueshovvuilter}. In particular, we established the existence of a global
attractor, determined its detailed structure and were able to compute the
Lyapunov exponents explicitly in some cases. We then extended these results
to the case of parabolic SPDEs driven by a homogeneous multiplicative white
noise defined in Stratonovitch's sense in \cite{chueshovvuilquinto},
investigated there various stability properties of the non-random global
attractor and established the existence of a recurrent motion of sorts among
its components. Furthermore, in \cite{chueshovvuilsexto} we analyzed the
same type of equations as in \cite{chueshovvuilquinto} but with the noise
defined in It\^{o}'s sense. In this way we were able to establish the
existence and many properties of a random global attractor\textit{\ }and
excluded in particular the existence of any kind of recurrence phenomena,
thereby obtaining radically different results than in \cite%
{chueshovvuilquinto}. The analysis carried out in \cite{chueshovvuilsexto}
was further deepened in \cite{bergechueshovvuilbis}, where it was shown that
the stabilization of the solution random fields toward the global attractor
is entirely controlled by their spatial average, thereby obtaining exchange
of stability results particularly relevant to the description of certain
migration phenomena in population dynamics. Finally, in \cite%
{chueshovvuilseptimo} we proved the existence of invariant sets under the
flow generated by certain systems of SPDEs including those of Lotka-Volterra
and Landau-Ginzburg.

But Igor's interests did not limit themselves to investigations of solutions
to SPDEs as he was also genuinely interested in the many possible
connections that exist between systems of differential equations on the one
hand, and the theory of random dynamical systems and stochastic processes on
the other hand (see, e.g., \cite{chueshov}). This prompted me to present
here some very recent and preliminary results concerning the time evolution
of certain Bernstein processes naturally associated with a class of
deterministic linear partial differential equations. Accordingly, the
remaining part of this article is organized as follows: In Section 2 I
recall what a Bernstein process is, and state there a theorem that shows how
to associate such a process with the two adjoint parabolic Cauchy problems%
\begin{eqnarray}
\partial _{t}u(\mathsf{x},t) &=&\frac{1}{2}\triangle _{\mathsf{x}}u(\mathsf{x%
},t)-V\left( \mathsf{x}\right) u(\mathsf{x},t),\text{ \ \ }(\mathsf{x},t)\in 
\mathbb{R}^{d}\mathbb{\times }\left( 0,T\right] ,  \notag \\
u(\mathsf{x},0) &=&\varphi (\mathsf{x})=\mathcal{N}\varphi _{0}(\mathsf{x}),%
\text{ \ }\mathsf{x}\in \mathbb{R}^{d}  \label{forwardcauchy1}
\end{eqnarray}%
and%
\begin{eqnarray}
-\partial _{t}v(\mathsf{x},t) &=&\frac{1}{2}\triangle _{\mathsf{x}}v(\mathsf{%
x},t)-V\left( \mathsf{x}\right) v(\mathsf{x},t),\text{ \ \ }(\mathsf{x}%
,t)\in \mathbb{R}^{d}\mathbb{\times }\left[ 0,T\right) ,  \notag \\
v(\mathsf{x},T) &=&\psi (\mathsf{x})=\mathcal{N}\psi _{T}(\mathsf{x}),\text{
\ }\mathsf{x}\in \mathbb{R}^{d},  \label{backwardcauchy1}
\end{eqnarray}%
where $T$ $>0$ is arbitrary and where $\triangle _{\mathsf{x}}$ stands for
Laplace's operator with respect to the spatial variable. In these equations $%
\mathcal{N}>0$ is a normalization factor whose significance I explain below.
Moreover, $V$ is real-valued while $\varphi _{0}$ and $\psi _{T}$ are
positive data which are assumed to be either Gaussian functions of the form%
\begin{eqnarray}
\varphi _{0}(\mathsf{x}) &=&\exp \left[ -\frac{\left\vert \mathsf{x-a}%
_{0}\right\vert ^{2}}{2\sigma _{0}}\right] ,  \label{gaussian1} \\
\psi _{T}(\mathsf{x}) &=&\exp \left[ -\frac{\left\vert \mathsf{x-a}%
_{T}\right\vert ^{2}}{2\sigma _{T}}\right]  \label{gaussian2}
\end{eqnarray}%
where $\sigma _{0,T}>0$ and $\mathsf{a}_{0,T}\in \mathbb{R}^{d}$ are
arbitrary vectors with $\left\vert .\right\vert $ the usual Euclidean norm,
or%
\begin{eqnarray}
\varphi _{0}(\mathsf{x}) &=&\dprod\limits_{j=1}^{d}\left( \left( 1-\frac{%
\left\vert x_{j}-a_{0,j}\right\vert }{\sigma _{0}}\right) \vee 0\right) ,
\label{hat1} \\
\psi _{T}(\mathsf{x}) &=&\dprod\limits_{j=1}^{d}\left( \left( 1-\frac{%
\left\vert x_{j}-a_{T,j}\right\vert }{\sigma _{T}}\right) \vee 0\right) .
\label{hat2}
\end{eqnarray}%
In (\ref{hat1})-(\ref{hat2}), $x_{j}$ and $a_{0,T,j}$ denote the $j^{th}$
component of $\mathsf{x}$ and $\mathsf{a}_{0,T}$, respectively. Furthermore
these initial-final conditions have localization properties which are more
clear-cut than those of (\ref{gaussian1})-(\ref{gaussian2}) in that they
vanish identically outside hypercubes in $\mathbb{R}^{d}$. The cases where

\begin{equation}
\varphi _{0}(\mathsf{x})=\delta _{0}(\mathsf{x})  \label{dirac}
\end{equation}%
with $\delta _{0}$ the Dirac measure concentrated at the origin and $\psi
_{T}$ given by (\ref{gaussian2}) or (\ref{dirac}) are also considered. An
important observation here is that (\ref{gaussian1})-(\ref{gaussian2}) and (%
\ref{hat1})-(\ref{hat2}) are not normalized as standard probability
distributions, for the only normalization condition needed below involves $%
\varphi _{0}$, $\psi _{T}$ and $\mathcal{N}$ in a rather unexpected way
which is inherently tied up with the construction of Bernstein processes.
Finally, the following hypothesis is imposed regarding the potential
function in (\ref{forwardcauchy1})-(\ref{backwardcauchy1}):

\bigskip

(H) The function $V:\mathbb{R}^{d}\mapsto \mathbb{R}$ is continuous, bounded
from below and satisfies $V\left( \mathsf{x}\right) \rightarrow +\infty $ as 
$\left\vert \mathsf{x}\right\vert \rightarrow +\infty $.

\bigskip

An immediate consequence of this hypothesis is that the resolvent of the
usual self-adjoint realization of the elliptic operator on the right-hand
side of (\ref{forwardcauchy1})-(\ref{backwardcauchy1}) is compact in $L_{%
\mathbb{C}}^{2}\left( \mathbb{R}^{d}\right) $, the usual Lebesgue space of
all square integrable, complex-valued functions on $\mathbb{R}^{d}$. This
means that the operator in question has an entirely discrete spectrum $%
\left( E_{\mathsf{n}}\right) _{\mathsf{n}\in \mathbb{N}^{d}}$, and that
there exists an orthonormal basis $\left( \mathsf{f}_{\mathsf{n}}\right) _{%
\mathsf{n}\in \mathbb{N}^{d}}\subset L_{\mathbb{C}}^{2}\left( \mathbb{R}%
^{d}\right) $ consisting entirely of its eigenfunctions (see, e.g., Section
XIII.14 in \cite{reedsimon}). In the context of this article the convergence
of the series%
\begin{equation}
\sum_{\mathsf{n}\in \mathbb{N}^{d}}\exp \left[ -tE_{\mathsf{n}}\right]
<+\infty  \label{convergence}
\end{equation}%
for every $t\in \left( 0,T\right] $ is also required. Then, under the above
conditions the construction of a Markovian Bernstein process rests on two
essential ingredients, namely, Green's function (or heat kernel) associated
with (\ref{forwardcauchy1})-(\ref{backwardcauchy1}), which satisfies the
symmetry and positivity conditions%
\begin{equation}
g(\mathsf{x},t,\mathsf{y})=g(\mathsf{y},t,\mathsf{x})>0
\label{greenfunction}
\end{equation}%
for all $\mathsf{x},\mathsf{y\in }\mathbb{\ \mathbb{R}}^{d}\mathbb{\ }$and
every $t\in \left( 0,T\right] $, and the probability measure $\mu $ on $%
\mathbb{\mathbb{R}}^{d}\times \mathbb{\mathbb{R}}^{d}$ whose density is
given by%
\begin{equation}
\mu (\mathsf{x,y})=\varphi (\mathsf{x)}g(\mathsf{x},T,\mathsf{y})\psi (%
\mathsf{y),\label{density}}
\end{equation}%
which satisfies the normalization condition%
\begin{eqnarray}
&&\int_{\mathbb{R}^{d}\times \mathbb{R}^{d}}\mathsf{dxdy}\varphi (\mathsf{x)}%
g(\mathsf{x},T,\mathsf{y})\psi (\mathsf{y)}  \notag \\
&=&\mathcal{N}^{2}\int_{\mathbb{R}^{d}\times \mathbb{R}^{d}}\mathsf{dxdy}%
\varphi _{0}(\mathsf{x)}g(\mathsf{x},T,\mathsf{y})\psi _{T}(\mathsf{y)=1.%
\label{normalization}}
\end{eqnarray}%
Notice that (\ref{normalization}) may be considered as the definition of $%
\mathcal{N}$, and that the inequality in (\ref{greenfunction}) is a
consequence of two-sided Gaussian bounds for $g$ whose existence follows
from the general theory developed in \cite{aronson} and further refined in
Chapter 3 of \cite{davies}. Moreover, as a consequence of (H) and (\ref%
{convergence}), Green's function admits an expansion of the form%
\begin{equation}
g(\mathsf{x},t,\mathsf{y})=\sum_{\mathsf{n}\in \mathbb{N}^{d}}\exp \left[
-tE_{\mathsf{n}}\right] \mathsf{f}_{\mathsf{n}}(\mathsf{x})\mathsf{f}_{%
\mathsf{n}}(\mathsf{y})  \label{expansion}
\end{equation}%
which converges strongly in $L_{\mathbb{C}}^{2}\left( \mathbb{R}^{d}\times 
\mathbb{R}^{d}\right) $ for every $t\in \left( 0,T\right] $ (unless more
detailed information about the $\mathsf{f}_{\mathsf{n}}$'s or
ultracontractive bounds become available, in which case the convergence can
be substantially improved, see, e.g., Chapter 2 in \cite{davies}). Thus, in
Section 2 the knowledge of $g$ and $\mu $ is used to state a theorem about
the existence of a probability space which supports a Markovian Bernstein
process $Z_{\tau \in \left[ 0,T\right] }$ whose state space is the entire
Euclidean space $\mathbb{\mathbb{R}}^{d}$, and which is characterized by its
finite-dimensional distributions, the joint distribution of $Z_{0}$ and $%
Z_{T}$ and the probability of finding $Z_{t}$ at any time $t\in \left[ 0,T%
\right] $ in a given region of space. In that section a very simple result
regarding the time evolution of $Z_{\tau \in \left[ 0,T\right] }$ is also
proved when considering (\ref{forwardcauchy1})-(\ref{backwardcauchy1}) with (%
\ref{hat1})-(\ref{hat2}). Section 3 is devoted to the analysis of the
function that determines the time evolution of the probability of finding $%
Z_{\tau \in \left[ 0,T\right] }$ in particular two-dimensional geometric
shapes that exhibit spherical symmetry in the case of the so-called harmonic
potential%
\begin{equation}
V(\mathsf{x})=\frac{\left\vert \mathsf{x}\right\vert ^{2}}{2},
\label{harmonicpot}
\end{equation}%
and for various combinations of the initial-final data given above. Finally,
a simple Faedo-Galerkin scheme is proposed whose ultimate goal is to allow
approximate computations of all the probability distributions involved.

\section{An existence result for a class of Bernstein processes in $\mathbb{R%
}^{d}$}

As a stochastic process a Bernstein process may be defined independently of
any reference to a system of partial differential equations, and there are
several equivalent ways to do so (see, e.g., \cite{jamison}). I shall
restrict myself to the following:

\bigskip

\textbf{Definition.} Let $d\in \mathbb{N}^{+}$ and $T\in \left( 0,+\infty
\right) $ be arbitrary. An $\mathbb{R}^{d}$-valued process $Z_{\tau \in %
\left[ 0,T\right] }$ defined on the complete probability space $\left(
\Omega ,\mathcal{F},\mathbb{P}\right) $ is called a Bernstein process if 
\begin{equation}
\mathbb{E}\left( f(Z_{r})\left\vert \mathcal{F}_{s}^{+}\vee \mathcal{F}%
_{t}^{-}\right. \right) =\mathbb{E}\left( f(Z_{r})\left\vert
Z_{s},Z_{t}\right. \right)  \label{condiexpect}
\end{equation}%
for every bounded Borel measurable function $f:\mathbb{R}^{d}\mapsto \mathbb{%
R}$ and for all $r,s,t$ satisfying $r\in \left( s,t\right) \subset \left[ 0,T%
\right] $. In (\ref{condiexpect}), the $\sigma $-algebras are%
\begin{equation}
\mathcal{F}_{s}^{+}=\sigma \left\{ Z_{\tau }^{-1}\left( F\right) :\tau \leq
s,\text{ }F\in \mathcal{B}_{d}\right\}  \label{pastalgebra}
\end{equation}%
and%
\begin{equation}
\mathcal{F}_{t}^{-}=\sigma \left\{ Z_{\tau }^{-1}\left( F\right) :\tau \geq
t,\text{ }F\in \mathcal{B}_{d}\right\} ,  \label{futurealgebra}
\end{equation}%
where $\mathcal{B}_{d}$ stands for the Borel $\sigma $-algebra on $\mathbb{R}%
^{d}$. Moreover, $\mathbb{E}$ denotes the (conditional) expectation
functional on $\left( \Omega ,\mathcal{F},\mathbb{P}\right) $.

\bigskip

The dynamics of such a process are, therefore, solely determined by the
properties of the process at times $s$ and $t$, irrespective of its behavior
prior to instant $s$ and after instant $t$. Of course, it is plain that this
fact generalizes the usual Markov property.

In what follows an important r\^{o}le is played by the positive solution to (%
\ref{forwardcauchy1}) and the positive solution to (\ref{backwardcauchy1}),
namely,%
\begin{equation}
u(\mathsf{x},t)=\int_{\mathbb{R}^{d}}\mathsf{dy}g(\mathsf{x},t,\mathsf{y}%
)\varphi (\mathsf{y)\label{forwardsolution}}
\end{equation}%
and 
\begin{equation}
v(\mathsf{x},t)=\int_{\mathbb{R}^{d}}\mathsf{dy}g(\mathsf{x},T-t,\mathsf{y}%
)\psi (\mathsf{y),\label{backwardsolution}}
\end{equation}%
respectively. Taken together, (\ref{forwardcauchy1}) and (\ref%
{backwardcauchy1}) may thus be looked upon as defining a decoupled
forward-backward system of linear deterministic partial differential
equations, with (\ref{forwardsolution}) wandering off to the future and (\ref%
{backwardsolution}) evolving into the past. The functions%
\begin{equation}
p\left( \mathsf{x},t;\mathsf{z},r;\mathsf{y},s\right) =g^{-1}(\mathsf{x},t-s,%
\mathsf{y})g(\mathsf{x},t-r,\mathsf{z})g(\mathsf{z},r-s,\mathsf{y})
\label{bernsteindensity}
\end{equation}%
and%
\begin{equation}
P\left( \mathsf{x},t;F,r;\mathsf{y},s\right) =\dint\limits_{F}\mathsf{dz}%
p\left( \mathsf{x},t;\mathsf{z},r;\mathsf{y},s\right)  \label{transition}
\end{equation}%
with $F\in \mathcal{B}_{d}$, both being well defined and positive for all $%
\mathsf{x},\mathsf{y},\mathsf{z}\in \mathbb{R}^{d}$ and all $r,s,t$
satisfying $r\in \left( s,t\right) \subset \left[ 0,T\right] $, are equally
important as is the probability measure $\mu $ whose density is (\ref%
{density}), namely,%
\begin{equation}
\mu \left( G\right) =\int_{G}\mathsf{dxdy}\varphi (\mathsf{x)}g(\mathsf{x},T,%
\mathsf{y})\psi (\mathsf{y)\label{probabilitymeasure}}
\end{equation}%
where $G\in \mathcal{B}_{d}\times \mathcal{B}_{d}$, which satisfies the
normalization condition (\ref{normalization}). The corresponding initial and
final marginal distributions then read%
\begin{eqnarray*}
\mu \left( F\times \mathbb{R}^{d}\right) &=&\int_{F}\mathsf{dx}\varphi (%
\mathsf{x})\int_{\mathbb{R}^{d}}\mathsf{dy}g(\mathsf{x},T,\mathsf{y})\psi (%
\mathsf{y}) \\
&=&\int_{F}\mathsf{dx}\varphi (\mathsf{x})v(\mathsf{x},0)
\end{eqnarray*}%
and%
\begin{eqnarray*}
\mu (\mathbb{R}^{d}\times F) &=&\int_{F}\mathsf{dy}\psi (\mathsf{y})\int_{%
\mathbb{R}^{d}}\mathsf{dx}g(\mathsf{x},T,\mathsf{y})\varphi (\mathsf{x}) \\
&=&\int_{F}\mathsf{dy}u(\mathsf{y},T)\psi (\mathsf{y})
\end{eqnarray*}%
respectively, as a consequence of (\ref{forwardsolution}) and (\ref%
{backwardsolution}). It is the knowledge of (\ref{transition}) and (\ref%
{probabilitymeasure}) that makes it possible to associate with (\ref%
{forwardcauchy1}) and (\ref{backwardcauchy1}) a Bernstein process in the
following sense:

\bigskip

\textbf{Theorem.} \textit{Assume that }$V$\textit{\ satisfies} \textit{%
Hypothesis (H), that condition (\ref{convergence}) holds and that }$P$%
\textit{\ and }$\mu $\textit{\ are given by (\ref{transition}) and (\ref%
{probabilitymeasure}), respectively. Then there exists a probability space }$%
\left( \Omega ,\mathcal{F},\mathbb{P}_{\mu }\right) $\textit{\ supporting an 
}$\mathbb{R}^{d}$\textit{-valued Bernstein process }$Z_{\tau \in \left[ 0,T%
\right] }$\textit{\ such that the following properties are valid:}

\textit{(a) The process }$Z_{\tau \in \left[ 0,T\right] }$\textit{\ is
Markovian, and the function }$P$\textit{\ is its transition function in the
sense that}%
\begin{equation*}
\mathbb{P}_{\mu }\left( Z_{r}\in F\left\vert Z_{s},Z_{t}\right. \right)
=P\left( Z_{t},t;F,r;Z_{s},s\right)
\end{equation*}%
\textit{for each }$F\in \mathcal{B}_{d}$\textit{\ and all }$r,s,t$\textit{\
satisfying }$r\in \left( s,t\right) \subset \left[ 0,T\right] $\textit{.
Moreover,}%
\begin{equation}
\mathbb{P}_{\mu }\left( Z_{0}\in F_{0},Z_{T}\in F_{T}\right) =\mu
(F_{0}\times F_{T})  \label{jointproba}
\end{equation}%
\textit{for all }$F_{0},F_{T}\in \mathcal{B}_{d}$\textit{, that is, }$\mu $%
\textit{\ is the joint probability distribution of }$Z_{0}$\textit{\ and }$%
Z_{T}$\textit{.}

\textit{(b) The finite-dimensional probability distributions of the process
are given by}%
\begin{eqnarray}
&&\mathbb{P}_{\mu }\left( Z_{t_{1}}\in F_{1},...,Z_{t_{n}}\in F_{n}\right)
\label{projections} \\
&=&\int_{F_{1}}\mathsf{dx}_{1}...\int_{F_{n}}\mathsf{dx}_{n}\dprod%
\limits_{k=2}^{n}g\left( \mathsf{x}_{k},t_{k}-t_{k-1},\mathsf{x}%
_{k-1}\right) \times u(\mathsf{x}_{1},t_{1})v(\mathsf{x}_{n},t_{n})  \notag
\end{eqnarray}%
\textit{for every integer }$n\geq 2$\textit{, all }$F_{1},...,F_{n}\in 
\mathcal{B}_{d}$\textit{\ and all }$t_{0}=0<t_{1}<...<t_{n}<T$\textit{,
where }$u$ \textit{and} $v$ \textit{are given by (\ref{forwardsolution}) and
(\ref{backwardsolution}), respectively.}

\textit{(c) The probability of finding the process in a given region }$%
F\subset \mathbb{R}^{d}$\textit{\ at time }$t$ \textit{is given by}%
\begin{equation}
\mathbb{P}_{\mu }\left( Z_{t}\in F\right) =\int_{F}\mathsf{dx}u(\mathsf{x}%
,t)v(\mathsf{x},t)  \label{probability}
\end{equation}%
\textit{for each }$F\in \mathcal{B}_{d}$\textit{\ and every }$t\in \left[ 0,T%
\right] .$

\textit{(d) }$\mathbb{P}_{\mu }$\textit{\ is the only probability measure
leading to the above properties.}

\bigskip

I omit the proof of this theorem, which can be adapted either from the
abstract arguments in \cite{jamison} or from the more analytical approach in 
\cite{vuillerzambrin1}, and will rather focus on its consequences regarding
the time evolution of $Z_{\tau \in \left[ 0,T\right] }$. Prior to that some
comments are in order:

\bigskip

\textsc{Remarks.} (1) Hypothesis (H) and condition (\ref{convergence}) are
sufficient but not necessary for the theorem to hold. However, the advantage
of having (\ref{expansion}) is that such an expansion greatly simplifies
some calculations and also has the virtue of making theoretical results
amenable to approximations and computations. I will dwell a bit more on this
point in the next section.

(2) Bernstein processes may be Markovian but in general they are not.
Independently of that they have played an increasingly important r\^{o}le in
various areas of mathematics and physics over the years. It is not possible
to give a complete bibliography here, but I will refer instead the
interested reader to \cite{jamison}, \cite{roellythieullen} and \cite%
{vuillerzambrin1} which contain many references describing the history and
earlier works on the subject, tracing things back to the pioneering works 
\cite{bernstein} and \cite{schroedinger}. Moreover, Bernstein processes have
also lurked in various forms in more recent applications of Optimal
Transport Theory, as testified by the monographs \cite{galichon} and \cite%
{villani}. In this regard it is worth mentioning that they are also referred
to as \textit{Schr\"{o}dinger processes} or \textit{reciprocal processes} in
the literature.

(3) The probability measure $\mu $ of a non-Markovian Bernstein process does
not have as simple a structure as that given by (\ref{probabilitymeasure}).
A case in point is the so-called periodic Ornstein-Uhlenbeck process, which
is one of the simplest stationary Gaussian non-Markovian processes that can
be viewed as a particular Bernstein process, as was recently proved in \cite%
{vuillerzambrin2} (see also, e.g., \cite{roellythieullen} and the references
therein for other analyses of the periodic Ornstein-Uhlenbeck process). In
this case the construction of the measure $\mu $ is much more complicated
than in the Markovian case, as it involves a weighted average of a sequence
of suitably constructed signed measures naturally associated with an
infinite hierarchy of forward-backward linear parabolic equations.

\bigskip

Coming back to the main theme of this article, it is interesting to note
that the probability of finding the process at any given time $t\in \left[
0,T\right] $ in an arbitrary region of space is expressed as an integral of
the product of $u$ and $v$ through the simple formula (\ref{probability}).
This is a manifestation of the fact that the process $Z_{\tau \in \left[ 0,T%
\right] }$ is actually reversible and exhibits a perfect symmetry between
past and future, a property already built to some extent into the definition
given at the beginning of this section. It is of course difficult to say
more about the time evolution of $Z_{\tau \in \left[ 0,T\right] }$ unless we
know more about the potential function $V$. However, at the very least the
following result holds, which in effect describes a recurrence property of
the process in a particular case:

\bigskip

\textbf{Proposition 1.} \textit{Let }$Z_{\tau \in \left[ 0,T\right] }$ 
\textit{be the Bernstein process associated with (\ref{forwardcauchy1})-(\ref%
{backwardcauchy1}) in the sense of the above theorem, where }$\varphi _{0}$ 
\textit{and }$\psi _{T}$\textit{\ are given by (\ref{hat1}) and (\ref{hat2}%
), respectively, and let}%
\begin{equation*}
\mathsf{C}_{\mathsf{a}_{0},\sigma _{0}}=\left\{ \mathsf{x}\in \mathbb{R}%
^{d}:\left\vert x_{j}-a_{0,j}\right\vert <\sigma _{0},\text{ \ }%
j=1,...,d\right\}
\end{equation*}%
\textit{be the hypercube outside which }$\varphi _{0}$\textit{\ vanishes
identically, that is, }$\varphi _{0}=0$\textit{\ on }$F_{\mathsf{a}%
_{0}},_{\sigma _{0}}=\mathbb{R}^{d}\setminus \mathsf{C}_{\mathsf{a}%
_{0},\sigma _{0}}$\textit{. Let }$\mathsf{C}_{\mathsf{a}_{T},\sigma _{T}}$%
\textit{\ be defined in a similar} \textit{way}. \textit{Then }%
\begin{equation*}
\mathbb{P}_{\mu }\left( Z_{0}\in \mathsf{C}_{\mathsf{a}_{0},\sigma
_{0}}\right) =1
\end{equation*}%
\textit{and}%
\begin{equation*}
\mathbb{P}_{\mu }\left( Z_{T}\in \mathsf{C}_{\mathsf{a}_{T},\sigma
_{T}}\right) =1.
\end{equation*}

\bigskip

\textbf{Proof.} This is an immediate consequence of (\ref{probability}), for%
\begin{equation*}
\mathbb{P}_{\mu }\left( Z_{0}\in F_{\mathsf{a}_{0}},_{\sigma _{0}}\right)
=\int_{F_{\mathsf{a}_{0}},_{\sigma _{0}}}\mathsf{dx}\varphi (\mathsf{x}%
)v\left( \mathsf{x},0\right) =0
\end{equation*}%
and%
\begin{equation*}
\mathbb{P}_{\mu }\left( Z_{T}\in F_{\mathsf{a}_{T}},_{\sigma _{T}}\right)
=\int_{F_{\mathsf{a}_{T}},_{\sigma _{T}}}\mathsf{dx}u\left( \mathsf{x}%
,T\right) \psi (\mathsf{x})=0\text{. \ \ }\blacksquare
\end{equation*}

\bigskip

Thus, in this case the process certainly starts its journey within $\mathsf{C%
}_{\mathsf{a}_{0},\sigma _{0}}$ and ends it within $\mathsf{C}_{\mathsf{a}%
_{T},\sigma _{T}}$. Since this is true no matter how small $\sigma _{0,T}$
are, that constitutes a generalization of the class of Bernstein bridges
constructed in \cite{vuillerzambrin2}. In particular, if $\mathsf{a}_{0}=%
\mathsf{a}_{T}$ and if $\sigma _{T}\leq \sigma _{0}$ the inclusion $\mathsf{C%
}_{\mathsf{a}_{T},\sigma _{T}}\subseteq \mathsf{C}_{\mathsf{a}_{0},\sigma
_{0}}$ holds, so that the process goes back to the region where it started
from with probability one, independently of its unknown whereabouts at
intermediary times $t\in \left( 0,T\right) $. These properties and
Proposition 1 remain true for all choices of $\varphi _{0}$, $\psi _{T}$
that vanish identically outside of a given Borel set, for instance for the
isotropic version of (\ref{hat1})-(\ref{hat2}), namely,%
\begin{eqnarray*}
\varphi _{0}(\mathsf{x}) &=&\left( 1-\frac{\left\vert \mathsf{x}-\mathsf{a}%
_{0}\right\vert }{\sigma _{0}}\right) \vee 0, \\
\psi _{T}(\mathsf{x}) &=&\left( 1-\frac{\left\vert \mathsf{x}-\mathsf{a}%
_{T}\right\vert }{\sigma _{T}}\right) \vee 0,
\end{eqnarray*}%
provided the sets $\mathsf{C}_{\mathsf{a}_{0,T},\sigma _{0,T}}$ are replaced
by the $d$-dimensional open balls%
\begin{equation*}
\mathsf{B}_{\mathsf{a}_{0,T},\sigma _{0,T}}=\left\{ \mathsf{x}\in \mathbb{R}%
^{d}:\left\vert \mathsf{x}-\mathsf{a}_{0,T}\right\vert <\sigma _{0,T}\right\}
\end{equation*}%
of radius $\sigma _{0,T}$ centered at $\mathsf{a}_{0,T}$. It would be
interesting to carry out a numerical simulation in real time of the behavior
of the processes generated in this way.

The preceding result fails to hold if the initial-final data are not of the
above form. In the next section I investigate this issue more closely in
case the potential function is given by (\ref{harmonicpot}).

\section{Some new results for the harmonic case}

The starting point is thus the forward-backward system%
\begin{eqnarray}
\partial _{t}u(\mathsf{x},t) &=&\frac{1}{2}\triangle _{\mathsf{x}}u(\mathsf{x%
},t)-\frac{\left\vert \mathsf{x}\right\vert ^{2}}{2}u(\mathsf{x},t),\text{ \
\ }(\mathsf{x},t)\in \mathbb{R}^{d}\mathbb{\times }\left( 0,T\right] , 
\notag \\
u(\mathsf{x},0) &=&\varphi \left( \mathsf{x}\right) \mathsf{=}\mathcal{N}%
\varphi _{0}\left( \mathsf{x}\right) ,\text{ \ }\mathsf{x}\in \mathbb{R}^{d}
\label{forwardcauchy2}
\end{eqnarray}%
and%
\begin{eqnarray}
-\partial _{t}v(\mathsf{x},t) &=&\frac{1}{2}\triangle _{\mathsf{x}}v(\mathsf{%
x},t)-\frac{\left\vert \mathsf{x}\right\vert ^{2}}{2}v(\mathsf{x},t),\text{
\ \ }(\mathsf{x},t)\in \mathbb{R}^{d}\mathbb{\times }\left[ 0,T\right) , 
\notag \\
v(\mathsf{x},T) &=&\psi \left( \mathsf{x}\right) \mathsf{=}\mathcal{N}\psi
_{T}\left( \mathsf{x}\right) ,\text{ \ }\mathsf{x}\in \mathbb{R}^{d}.
\label{backwardcauchy2}
\end{eqnarray}%
Green's function associated with (\ref{forwardcauchy2})-(\ref%
{backwardcauchy2}) is known to be Mehler's multidimensional kernel%
\begin{eqnarray}
&&g(\mathsf{x},t,\mathsf{y})  \notag \\
&=&\left( 2\pi \sinh \left( t\right) \right) ^{-\frac{d}{2}}\exp \left[ -%
\frac{\cosh \left( t\right) \left( \left\vert \mathsf{x}\right\vert
^{2}+\left\vert \mathsf{y}\right\vert ^{2}\right) -2\left( \mathsf{x,y}%
\right) _{\mathbb{R}^{d}}}{2\sinh \left( t\right) }\right]  \label{mehler}
\end{eqnarray}%
where $(.,.)_{\mathbb{R}^{d}}$ denotes the usual inner product in $\mathbb{%
\mathbb{R}}^{d}$ (see, e.g., the Appendix in \cite{vuillerzambrin2}). Then
if $\varphi _{0}$, $\psi _{T}$ are given by (\ref{gaussian1})-(\ref%
{gaussian2}), the solutions (\ref{forwardsolution})-(\ref{backwardsolution})
and the integral on the left-hand side of (\ref{normalization}) can all be
computed explicitly since the integrals are Gaussian. For instance, the
forward solution reads%
\begin{eqnarray}
&&u(\mathsf{x,}t)=\mathcal{N}\left( \frac{\sigma _{0}}{\sigma _{0}\cosh
(t)+\sinh (t)}\right) ^{\frac{d}{2}}\exp \left[ -\frac{\left\vert \mathsf{a}%
_{0}\right\vert ^{2}}{2\sigma _{0}}\right]  \notag \\
&&\times \exp \left[ -\frac{\cosh (t)\left\vert \mathsf{x}\right\vert ^{2}}{%
2\sinh (t)}+\frac{\left\vert \sigma _{0}\mathsf{x+}\sinh (t)\mathsf{a}%
_{0}\right\vert ^{2}}{2\sigma _{0}\sinh (t)\left( \sigma _{0}\cosh (t)+\sinh
(t)\right) }\right]  \label{forwardsolutionbis}
\end{eqnarray}%
for every $t\in \left( 0,T\right] $, while the backward solution is obtained
from (\ref{forwardsolutionbis}) by replacing $\sigma _{0}$ by $\sigma _{T}$, 
$\mathsf{a}_{0}$ by $\mathsf{a}_{T}$ and $t$ by $T-t$, respectively. The
downside is that these expressions are complicated, cumbersome and in any
case unsuited to extract valuable information out of (\ref{probability})
unless particular choices are made for these parameters. For example, if $%
\sigma _{0}=\sigma _{T}=1$ and $\mathsf{a}_{0}=\mathsf{a}_{T}=0$, the
forward solution (\ref{forwardsolutionbis}) and the related backward
solution reduce to%
\begin{eqnarray}
&&u(\mathsf{x,}t)=\mathcal{N}\exp \left[ -\frac{\left\vert \mathsf{x}%
\right\vert ^{2}+dt}{2}\right] ,  \label{forwardsolutionter} \\
&&v(\mathsf{x,}t)=\mathcal{N}\exp \left[ -\frac{\left\vert \mathsf{x}%
\right\vert ^{2}+d(T-t)}{2}\right] ,  \label{backwardsolutionter}
\end{eqnarray}%
respectively, while an explicit computation from (\ref{normalization}) gives 
\begin{equation*}
\mathcal{N=\pi }^{-\frac{d}{4}}\exp \left[ \frac{dT}{4}\right]
\end{equation*}%
for the corresponding normalization factor. Therefore, the substitution of
these expressions into (\ref{probability}) leads to%
\begin{equation*}
\mathbb{P}_{\mu }\left( Z_{t}\in F\right) =\int_{F}\mathsf{dx}u(\mathsf{x}%
,t)v(\mathsf{x,}t)=\pi ^{-\frac{d}{2}}\int_{F}\mathsf{dx}\exp \left[
-\left\vert \mathsf{x}\right\vert ^{2}\right]
\end{equation*}%
for each $t\in \left[ 0,T\right] $ and every $F\in \mathcal{B}_{d}$, so that
the probability of finding the process in any region of space is here
independent of time. The reason for this independence can easily be
understood by means of the substitution of (\ref{mehler}) and (\ref%
{forwardsolutionter})-(\ref{backwardsolutionter}) into (\ref{projections}),
which first leads to the Gaussian law of $\left(
Z_{t_{1}},...,Z_{t_{n}}\right) \in \mathbb{R}^{nd}$ and from there
eventually to the covariance%
\begin{equation*}
\mathbb{E}_{\mu }\left( Z_{s}^{i}Z_{t}^{j}\right) =\frac{1}{2}\exp \left[
-\left\vert t-s\right\vert \right] \delta _{i,j}
\end{equation*}%
for all $s,t\in \left[ 0,T\right] $ and all $i,j\in \left\{ 1,...,d\right\} $%
, where $\mathbb{E}_{\mu }$ denotes the expectation functional on the
probability space of the theorem. Therefore, the Bernstein process thus
constructed identifies in law with the standard $d$-dimensional
Ornstein-Uhlenbeck velocity process, so that the choice of (\ref{gaussian1}%
)-(\ref{gaussian2}) as initial-final data corresponds in a sense to an
equilibrium situation whereby the law remains stationary (see, e.g., \cite%
{karatzasshreve} for general properties of this and related processes). For
instance, if%
\begin{equation*}
A_{R_{1},R_{2}}=\left\{ \mathsf{x\in }\mathbb{R}^{2}:R_{1}\leq \left\vert 
\mathsf{x}\right\vert <R_{2}\right\}
\end{equation*}%
is the two-dimensional annulus centered at the origin with $R_{1}\geq 0$ and 
$R_{2}>0$, then%
\begin{equation*}
\mathbb{P}_{\mu }\left( Z_{t}\in A_{R_{1},R_{2}}\right) =\exp \left[
-R_{1}^{2}\right] -\exp \left[ -R_{2}^{2}\right] .
\end{equation*}

The situation is quite different if the system (\ref{forwardcauchy2})-(\ref%
{backwardcauchy2}) is considered with $\varphi _{0}$ given by (\ref{dirac})
and $\psi _{T}$ given by (\ref{gaussian2}) where $\sigma _{T}=1$ and $%
a_{T}=0 $. In this case%
\begin{equation}
u(\mathsf{x,}t)=\mathcal{N}\left( 2\pi \sinh (t)\right) ^{-\frac{d}{2}}\exp %
\left[ -\frac{\coth (t)\left\vert \mathsf{x}\right\vert ^{2}}{2}\right]
\label{forwardsolutionquarto}
\end{equation}%
and%
\begin{equation}
v(\mathsf{x,}t)=\mathcal{N}\exp \left[ -\frac{\left\vert \mathsf{x}%
\right\vert ^{2}+d(T-t)}{2}\right]  \label{backwardsolutionquarto}
\end{equation}%
for the forward and backward solutions, respectively, and furthermore the
value of $\mathcal{N}$ can again be determined directly from (\ref%
{normalization}). Indeed the relevant integral is%
\begin{equation*}
\int_{\mathbb{R}^{d}\times \mathbb{R}^{d}}\mathsf{dxdy}\delta _{0}(\mathsf{x)%
}g(\mathsf{x},T,\mathsf{y})\exp \left[ -\frac{\left\vert \mathsf{y}%
\right\vert ^{2}}{2}\right] =\exp \left[ -\frac{dT}{2}\right]
\end{equation*}%
by virtue of (\ref{mehler}), so that%
\begin{equation*}
\mathcal{N}=\exp \left[ \frac{dT}{4}\right] .
\end{equation*}%
Therefore, one obtains in particular%
\begin{equation*}
\mathbb{P}_{\mu }\left( Z_{0}\in \mathbb{R}^{d}\backslash \left\{ \mathsf{o}%
\right\} \right) =\mathcal{N}\int_{\mathbb{R}^{d}\backslash \left\{ \mathsf{o%
}\right\} }\mathsf{dx}\delta _{0}(\mathsf{x)}v(\mathsf{x,0})=0
\end{equation*}%
so that the process is conditioned to start at the origin since%
\begin{equation}
\mathbb{P}_{\mu }\left( Z_{0}=\mathsf{o}\right) =1.  \label{pinneddown}
\end{equation}%
Moreover, for positive times an explicit evaluation from (\ref{probability})
leads to%
\begin{equation*}
\mathbb{P}_{\mu }\left( Z_{t}\in F\right) =\left( 2\pi \rho (t)\right) ^{-%
\frac{d}{2}}\int_{F}\mathsf{dx}\exp \left[ -\frac{\left\vert \mathsf{x}%
\right\vert ^{2}}{2\rho (t)}\right]
\end{equation*}%
where the width parameter is identified as%
\begin{equation}
\rho (t)=\sinh (t)\exp \left[ -t\right] .  \label{variance}
\end{equation}%
It is then instructive to consider again the case of $Z_{\tau \in \ \left[
0,T\right] }$ wandering in the two-dimensional annulus $A_{R_{1},R_{2}}$,
and to investigate the way that%
\begin{equation}
\mathbb{P}_{\mu }\left( Z_{t}\in A_{R_{1},R_{2}}\right) =\exp \left[ -\frac{%
R_{1}^{2}}{2\rho (t)}\right] -\exp \left[ -\frac{R_{2}^{2}}{2\rho (t)}\right]
\label{probabilityquarto}
\end{equation}%
varies in the course of time for various values of the radii:

\bigskip

\textbf{Proposition 2.} \textit{The following statements hold:}

\textit{(a) If }$0=R_{1}<R_{2}$ \textit{one has}%
\begin{equation*}
\mathbb{P}_{\mu }\left( Z_{0}\in A_{0,R_{2}}\right) =1
\end{equation*}%
\textit{and the function} $t\mapsto $\textit{\ }$\mathbb{P}_{\mu }\left(
Z_{t}\in A_{0,R_{2}}\right) $\textit{\ is monotone decreasing on }$\left[ 0,T%
\right] $, \textit{eventually} \textit{reaching the minimal value}%
\begin{equation*}
\mathbb{P}_{\mu }\left( Z_{T}\in A_{0,R_{2}}\right) =1-\exp \left[ -\frac{%
R_{2}^{2}}{2\rho (T)}\right] .
\end{equation*}

\textit{(b) If }$0<R_{1}<R_{2}<1$ \textit{one has}%
\begin{equation}
\mathbb{P}_{\mu }\left( Z_{0}\in A_{R_{1},R_{2}}\right) =0
\label{probabilityseptimo}
\end{equation}%
\textit{\ and}%
\begin{equation*}
\mathbb{P}_{\mu }\left( Z_{t}\in A_{R_{1},R_{2}}\right) >0
\end{equation*}%
\textit{as soon as} $t>0$. \textit{Moreover,} \textit{if }$T$\textit{\ is
sufficiently large there exists a }$t^{\ast }\in \left( 0,T\right) $\textit{%
\ such that the function }$t\mapsto $\textit{\ }$\mathbb{P}_{\mu }\left(
Z_{t}\in A_{R_{1},R_{2}}\right) $\textit{\ is monotone decreasing for every }%
$t\in \left[ t^{\ast },T\right] .$

\textit{(c) If }$1\leq R_{1}<R_{2}$ \textit{one still has (\ref%
{probabilityseptimo}), but the function} $t\mapsto $\textit{\ }$\mathbb{P}%
_{\mu }\left( Z_{t}\in A_{R_{1},R_{2}}\right) $\textit{\ is monotone
increasing throughout }$\left[ 0,T\right] $.

\bigskip

\textbf{Proof.} Statement (a) follows immediately from (\ref{pinneddown})
and (\ref{probabilityquarto}) for $R_{1}=0$, as does the very first part of
(b) since then $\mathsf{o}\notin A_{R_{1},R_{2}}$. Now%
\begin{equation*}
\frac{d}{dt}\mathbb{P}_{\mu }\left( Z_{t}\in A_{R_{1},R_{2}}\right) =\frac{%
\rho ^{\prime }(t)}{2\rho ^{2}(t)}\left( \chi (R_{1},t)-\chi (R_{2},t)\right)
\end{equation*}%
where%
\begin{equation}
\chi (R,t)=R^{2}\exp \left[ -\frac{R^{2}}{2\rho (t)}\right] ,
\label{auxfunction}
\end{equation}%
and for any fixed $t\in \left( 0,T\right] $ this function is monotone
increasing for $R<\sqrt{2\rho (t)}$ and monotone decreasing for $R>\sqrt{%
2\rho (t)}$. Furthermore, (\ref{variance}) and $t\mapsto $ $\sqrt{2\rho (t)}$
are monotone increasing and concave with $\sqrt{2\rho (t)}<1$ uniformly in $%
t $. Therefore, if $0<R_{1}<R_{2}<1$ and if $T$ is large enough, there
exists a $t^{\ast }\in \left( 0,T\right) $ such that $R_{1}<R_{2}<\sqrt{%
2\rho (t^{\ast })}\leq \sqrt{2\rho (t)}$ for every $t\in \left[ t^{\ast },T%
\right] $, which implies the last claim of (b) since then $\chi
(R_{1},t)-\chi (R_{2},t)<0$. Finally, if $1\leq R_{1}<R_{2}$ one has\textit{%
\ a fortiori }$\sqrt{2\rho (t)}<R_{1}<R_{2}$ for every $t\in \left[ 0,T%
\right] $ so that $\chi (R_{1},t)-\chi (R_{2},t)>0$, which implies (c). \ \ $%
\blacksquare $

\bigskip

A natural interpretation of Statement (a) is that the process leaves the
origin as soon as $t>0$, and tends to quickly "leak out" of the disk $%
A_{0,R_{2}}$ when $R_{2}$ is sufficiently small. Moreover, Statement (b)
means that the probability of finding the process in the annulus increases
for small times, then reaches a maximal value and eventually decreases for
large times when $R_{1}$ and $R_{2}$ are sufficiently small, in sharp
contrast to Statement (c) where the probability in question is monotone
increasing for all times if $R_{1}$ and $R_{2}$ are sufficiently large.
Finally, the substitution of (\ref{mehler}) and (\ref{forwardsolutionquarto}%
)-(\ref{backwardsolutionquarto}) into (\ref{projections}) again determines
the projection of the law onto $\mathbb{R}^{nd}$ and, after long algebraic
manipulations, the covariance%
\begin{equation*}
\mathbb{E}_{\mu }\left( Z_{s}^{i}Z_{t}^{j}\right) =\frac{1}{2}\exp \left[
-(t+s)\right] \left( \exp \left[ 2(t\wedge s)\right] -1\right) \delta _{i,j}
\end{equation*}%
for all $s,t\in \left[ 0,T\right] $ and all $i,j\in \left\{ 1,...,d\right\} $%
. Therefore, the Bernstein process thus constructed is identical in law with
the Ornstein-Uhlenbeck process conditioned to start at the origin of $%
\mathbb{R}^{d}$.

A last example can be provided by choosing $\varphi _{0}$ and $\psi _{T}$
both of the form (\ref{dirac}) in (\ref{forwardcauchy2})-(\ref%
{backwardcauchy2}). In this case one gets%
\begin{equation*}
u(\mathsf{x,}t)=\mathcal{N}\left( 2\pi \sinh (t)\right) ^{-\frac{d}{2}}\exp %
\left[ -\frac{\coth (t)\left\vert \mathsf{x}\right\vert ^{2}}{2}\right]
\end{equation*}%
and%
\begin{equation*}
v(\mathsf{x,}t)=\mathcal{N}\left( 2\pi \sinh (T-t)\right) ^{-\frac{d}{2}%
}\exp \left[ -\frac{\coth (T-t)\left\vert \mathsf{x}\right\vert ^{2}}{2}%
\right]
\end{equation*}%
for the respective solutions, where the exact value of the normalization
factor is%
\begin{equation*}
\mathcal{N=}\left( 2\pi \sinh (T)\right) ^{\frac{d}{4}}\text{.}
\end{equation*}%
Arguing as in the preceding example one then obtains%
\begin{equation}
\mathbb{P}_{\mu }\left( Z_{0}=\mathsf{o}\right) =\mathbb{P}_{\mu }\left(
Z_{T}=\mathsf{o}\right) =1  \label{pinneddownbis}
\end{equation}%
so that the process is conditioned to start and end at the origin, thereby
representing a random loop in $\mathbb{R}^{d}$. Moreover, for positive times
one still gets from (\ref{probability})%
\begin{equation*}
\mathbb{P}_{\mu }\left( Z_{t}\in F\right) =\left( 2\pi \rho (t)\right) ^{-%
\frac{d}{2}}\int_{F}\mathsf{dx}\exp \left[ -\frac{\left\vert \mathsf{x}%
\right\vert ^{2}}{2\rho (t)}\right]
\end{equation*}%
and in particular%
\begin{equation}
\mathbb{P}_{\mu }\left( Z_{t}\in A_{R_{1},R_{2}}\right) =\exp \left[ -\frac{%
R_{1}^{2}}{2\rho (t)}\right] -\exp \left[ -\frac{R_{2}^{2}}{2\rho (t)}\right]
\label{probabilityocto}
\end{equation}%
in the case of the two-dimensional annulus, but with a width parameter now
given by%
\begin{equation}
\rho (t)=\frac{\sinh (t)\sinh (T-t)}{\sinh (T)}  \label{variancebis}
\end{equation}%
for every $t\in \left[ 0,T\right] $. This function is quite different from (%
\ref{variance}), and the following result is valid:

\bigskip

\textbf{Proposition 3.} \textit{The following statements hold:}

\textit{(a) If }$0=R_{1}<R_{2}$ \textit{one has}%
\begin{equation}
\mathbb{P}_{\mu }\left( Z_{0}\in A_{0,R_{2}}\right) =\mathbb{P}_{\mu }\left(
Z_{T}\in A_{0,R_{2}}\right) =1.  \label{pinneddownter}
\end{equation}%
\textit{Moreover, the function} $t\mapsto $\textit{\ }$\mathbb{P}_{\mu
}\left( Z_{t}\in A_{0,R_{2}}\right) $\textit{\ is monotone decreasing on }$%
\left[ 0,\frac{T}{2}\right] $ \textit{and} \textit{monotone increasing on} $%
\left[ \frac{T}{2},T\right] $, \textit{thereby taking the minimal value}%
\begin{equation*}
\mathbb{P}_{\mu }\left( Z_{\frac{T}{2}}\in A_{0,R_{2}}\right) =1-\exp \left[
-\frac{R_{2}^{2}}{2\rho \left( \frac{T}{2}\right) }\right] .
\end{equation*}

\textit{(b) If }$1\leq R_{1}<R_{2}$ \textit{one has}%
\begin{equation*}
\mathbb{P}_{\mu }\left( Z_{0}\in A_{0,R_{2}}\right) =\mathbb{P}_{\mu }\left(
Z_{T}\in A_{0,R_{2}}\right) =0.
\end{equation*}%
\textit{Moreover, the function} $t\mapsto $\textit{\ }$\mathbb{P}_{\mu
}\left( Z_{t}\in A_{R_{1},R_{2}}\right) $\textit{\ is monotone increasing on 
}$\left[ 0,\frac{T}{2}\right] $ \textit{and} \textit{monotone decreasing on} 
$\left[ \frac{T}{2},T\right] $, \textit{thereby taking the maximal value}%
\begin{equation*}
\mathbb{P}_{\mu }\left( Z_{\frac{T}{2}}\in A_{R_{1},R_{2}}\right) =\exp %
\left[ -\frac{R_{1}^{2}}{2\rho \left( \frac{T}{2}\right) }\right] -\exp %
\left[ -\frac{R_{2}^{2}}{2\rho \left( \frac{T}{2}\right) }\right] .
\end{equation*}

\bigskip

\textbf{Proof. }While (\ref{pinneddownter}) follows from (\ref{pinneddownbis}%
), Relation (\ref{probabilityocto}) with $R_{1}=0$ leads to%
\begin{equation*}
\frac{d}{dt}\mathbb{P}_{\mu }\left( Z_{t}\in A_{0,R_{2}}\right) =-\frac{\rho
^{\prime }(t)}{2\rho ^{2}(t)}\chi (R_{2},t)
\end{equation*}%
where $\rho ^{\prime }(t)\geq 0$ for $t\in \left[ 0,\frac{T}{2}\right] $ and 
$\rho ^{\prime }(t)\leq 0$ for $t\in \left[ \frac{T}{2},T\right] $ according
to (\ref{variancebis}), which implies Statement (a). Statement (b) follows
from these properties of $\rho ^{\prime }$ and an analysis similar to that
of Statement (c) in Proposition 2. Indeed, we remark that the curve $\rho :%
\left[ 0,T\right] \mapsto \left[ 0,+\infty \right) $ given by (\ref%
{variancebis}) is concave aside from satisfying $\rho (0)=\rho (T)=0$, and
that it takes on the maximal value%
\begin{equation*}
\rho \left( \frac{T}{2}\right) =\frac{\sinh ^{2}\left( \frac{T}{2}\right) }{%
\sinh (T)}
\end{equation*}%
at the mid-point of the time interval. Therefore, the inequalities%
\begin{equation*}
\sqrt{2\rho (t)}\leq \sqrt{2\rho \left( \frac{T}{2}\right) }\leq 1
\end{equation*}%
hold for every $t\in \left[ 0,T\right] $, which implies that (\ref%
{auxfunction}) is monotone decreasing throughout the time interval as a
function of $R$, a consequence of the hypothesis regarding the radii. \ \ $%
\blacksquare $

\bigskip

The above properties of (\ref{variancebis}) thus show that the Bernstein
process of Proposition 3 constitutes a generalization of a Brownian loop,
that is, of a particular case of a Brownian bridge (see, e.g., \cite%
{karatzasshreve}). This renders the preceding result quite natural, in that
the probability of finding the process in the disk $A_{0,R_{2}}$ is minimal
at the mid-point of the time interval where there is maximal randomness. At
the same time, the situation is reversed if the annulus is relatively far
away from the origin.

As long as the regions of interest are spherically symmetric, the preceding
calculations may be performed in any dimension and not merely for $d=2$.
However, I shall refrain from doing that and rather focus briefly on what to
do when the values of the parameters $\sigma _{0,T}$ and $\mathsf{a}_{0,T}$
are arbitrary, or when other combinations of the above initial-final data
are chosen. It is here that an expansion of the form (\ref{expansion}) is
essential, and I will now show what (\ref{expansion}) reduces to in the case
of (\ref{mehler}). First, the spectral decomposition of the elliptic
operator on the right-hand side of (\ref{forwardcauchy2})-(\ref%
{backwardcauchy2}) is known explicitly (the operator identifies up to a sign
with the Hamiltonian of an isotropic system of quantum harmonic oscillators,
see, e.g.,\cite{messiah}). Indeed, let $\left( h_{n}\right) _{n\in \mathbb{N}%
}$ be the usual Hermite functions%
\begin{equation}
h_{n}(x)=\left( \pi ^{\frac{1}{2}}2^{n}n!\right) ^{-\frac{1}{2}}\exp \left[ -%
\frac{x^{2}}{2}\right] H_{n}(x)  \label{hermitefunctions}
\end{equation}%
where the $H_{n}$'s stand for the Hermite polynomials

\begin{equation}
H_{n}(x)=\left( -1\right) ^{n}\exp \left[ x^{2}\right] \frac{\mathsf{d}^{n}}{%
\mathsf{d}x^{n}}\exp \left[ -x^{2}\right] .  \label{hermitepolynomials}
\end{equation}%
Then, it is easily verified that the tensor products $\otimes
_{j=1}^{d}h_{n_{j}}$ where the $n_{j}$'s run independently over $\mathbb{N}$
provide an orthonormal basis of eigenfunctions in $L_{\mathbb{C}}^{2}\left( 
\mathbb{R}^{d}\right) $ which satisfy the eigenvalue equation%
\begin{equation*}
\left( -\frac{1}{2}\Delta _{\mathsf{x}}+\frac{\left\vert \mathsf{x}%
\right\vert ^{2}}{2}\right) \mathsf{h}_{\mathsf{n}}\left( \mathsf{x}\right)
=E_{\mathsf{n}}\mathsf{h}_{\mathsf{n}}\left( \mathsf{x}\right)
\end{equation*}%
for each $\mathsf{n\in }\mathbb{N}$ and every $\mathsf{x}\in \mathbb{R}^{d}$%
, where $\mathsf{n}=(n_{1},...,n_{d})\in \mathbb{N}^{d}$ and

\begin{eqnarray}
E_{\mathsf{n}} &=&\sum_{j=1}^{d}n_{j}+\frac{d}{2},  \label{energy} \\
\mathsf{h}_{\mathsf{n}} &=&\otimes _{j=1}^{d}h_{n_{j}}.
\label{tensorproduct}
\end{eqnarray}%
The immediate consequences are that (\ref{convergence}) holds, and that
expansion (\ref{expansion}) for (\ref{mehler}) takes the form%
\begin{equation}
g(\mathsf{x},t,\mathsf{y})=\sum_{\mathsf{n}\in \mathbb{N}^{d}}\exp \left[
-tE_{\mathsf{n}}\right] \mathsf{h}_{\mathsf{n}}\left( \mathsf{x}\right) 
\mathsf{h}_{\mathsf{n}}\left( \mathsf{y}\right)  \label{mehlerbis}
\end{equation}%
where the series is now absolutely convergent for each $t\in \left( 0,T%
\right] $ uniformly in all $\mathsf{x,y}\in \mathbb{R}^{d}$. This very last
statement follows from Cram\'{e}r-Charlier's inequality%
\begin{equation}
\left\vert \mathsf{h}_{\mathsf{n}}\left( \mathsf{x}\right) \mathsf{h}_{%
\mathsf{n}}\left( \mathsf{y}\right) \right\vert \leq k^{2d}\pi ^{-\frac{d}{2}%
}  \label{cramercharlier}
\end{equation}%
valid uniformly in $\mathsf{n}$, $\mathsf{x}$ and $\mathsf{y}$, where $k\leq
1.086435$ (see, e.g., Section 10.18 in \cite{emotricomi} and the references
therein).

The advantage of having (\ref{mehlerbis}) is that the forward solution (\ref%
{forwardsolution}) may now be rewritten in terms of the Fourier coefficients
of $\varphi $ and $\psi $ along the basis $\left( \mathsf{h}_{\mathsf{n}%
}\right) _{\mathsf{n}\in \mathbb{N}^{d}}$, namely,%
\begin{equation}
u(\mathsf{x,}t)=\sum_{\mathsf{n}\in \mathbb{N}^{d}}\alpha _{\mathsf{n}}\exp %
\left[ -tE_{\mathsf{n}}\right] \mathsf{h}_{\mathsf{n}}\left( \mathsf{x}%
\right)  \label{forwardfourier}
\end{equation}%
where%
\begin{equation}
\alpha _{\mathsf{n}}=\mathcal{N}\int_{\mathbb{R}^{d}}\mathsf{dx}\varphi _{0}(%
\mathsf{x})\mathsf{h}_{\mathsf{n}}\left( \mathsf{x}\right) ,
\label{fouriercoeff1}
\end{equation}%
which in case of Gaussian initial-final data provides a nice representation
of (\ref{forwardsolutionbis}). In a similar way the backward solution (\ref%
{backwardsolution}) is%
\begin{equation}
v(\mathsf{x,}t)=\sum_{\mathsf{n}\in \mathbb{N}^{d}}\beta _{\mathsf{n}}\exp %
\left[ -(T-t)E_{\mathsf{n}}\right] \mathsf{h}_{\mathsf{n}}\left( \mathsf{x}%
\right)  \label{backwardfourier}
\end{equation}%
where 
\begin{equation}
\beta _{\mathsf{n}}=\mathcal{N}\int_{\mathbb{R}^{d}}\mathsf{dx}\psi _{T}(%
\mathsf{x})\mathsf{h}_{\mathsf{n}}\left( \mathsf{x}\right) ,
\label{fouriercoeff2}
\end{equation}%
so that the normalization condition (\ref{normalization}) now reads%
\begin{equation}
\sum_{\mathsf{n}\in \mathbb{N}^{d}}\alpha _{\mathsf{n}}\exp \left[ -TE_{%
\mathsf{n}}\right] \beta _{\mathsf{n}}=1.  \label{normalizationbis}
\end{equation}%
This way of formulating things, in turn, leads to the possibility of
constructing a sequence of Faedo-Galerkin approximations to the problem at
hand. Thus for any positive integer $N\geq 1$, let $\mathsf{E}_{N}\left( 
\mathbb{R}^{d}\right) $ be the $N^{d}$-dimensional subspace of $L_{\mathbb{C}%
}^{2}\left( \mathbb{R}^{d}\right) $ generated by the $\mathsf{h}_{\mathsf{n}%
} $'s where $n_{j}\in \left\{ 0,...,N-1\right\} $ for each component of $%
\mathsf{n}$. Green's function (\ref{mehlerbis}) may then be approximated by%
\begin{equation}
g_{N}(\mathsf{x},t,\mathsf{y})=\sum_{\mathsf{n}:0\leq n_{j}\leq N-1}\exp %
\left[ -tE_{\mathsf{n}}\right] \mathsf{h}_{\mathsf{n}}\left( \mathsf{x}%
\right) \mathsf{h}_{\mathsf{n}}\left( \mathsf{y}\right)
\label{galerkinapprox1}
\end{equation}%
in $\mathsf{E}_{N}\left( \mathbb{R}^{d}\right) \otimes \mathsf{E}_{N}\left( 
\mathbb{R}^{d}\right) $, which leads to the approximations%
\begin{equation}
u_{N}(\mathsf{x,}t)=\sum_{\mathsf{n}:0\leq n_{j}\leq N-1}\alpha _{\mathsf{n}%
}\exp \left[ -tE_{\mathsf{n}}\right] \mathsf{h}_{\mathsf{n}}\left( \mathsf{x}%
\right)  \label{galerkinapprox2}
\end{equation}%
and%
\begin{equation}
v_{N}(\mathsf{x,}t)=\sum_{\mathsf{n}:0\leq n_{j}\leq N-1}\beta _{\mathsf{n}%
}\exp \left[ -(T-t)E_{\mathsf{n}}\right] \mathsf{h}_{\mathsf{n}}\left( 
\mathsf{x}\right)  \label{galerkinapprox3}
\end{equation}%
to (\ref{forwardfourier}) and (\ref{backwardfourier}), respectively.
Consequently, various numerical computations and controlled approximations
of the probability distributions of interest now become possible. I complete
this short article by a simple illustration of this fact stated in
Proposition 4 below, whose proof is based on the following result which
provides an approximate value for $\mathcal{N}$:

\bigskip

\textbf{Lemma.} \textit{Let (\ref{normalizationbis}) be written as}%
\begin{equation*}
\mathcal{N}^{2}\hat{\alpha}_{\mathsf{0}}\exp \left[ -TE_{\mathsf{0}}\right] 
\hat{\beta}_{\mathsf{0}}+\mathcal{N}^{2}\sum_{\mathsf{n}\in \mathbb{N}^{d},%
\text{ }\mathsf{n\neq 0}}\hat{\alpha}_{\mathsf{n}}\exp \left[ -TE_{\mathsf{n}%
}\right] \hat{\beta}_{\mathsf{n}}=1
\end{equation*}%
\textit{where}%
\begin{eqnarray*}
\hat{\alpha}_{\mathsf{n}} &=&\mathcal{N}^{-1}\alpha _{\mathsf{n}}, \\
\hat{\beta}_{\mathsf{n}} &=&\mathcal{N}^{-1}\beta _{\mathsf{n}}
\end{eqnarray*}%
\textit{for every }$\mathsf{n}\in \mathbb{N}^{d}$.\textit{\ Then for all }$%
\sigma _{0,T}>0$\textit{, }$a_{0,T}\in \mathbb{R}^{d}$\textit{, the unique
positive solution to }%
\begin{equation}
\mathcal{N}^{2}\hat{\alpha}_{\mathsf{0}}\exp \left[ -TE_{\mathsf{0}}\right] 
\hat{\beta}_{\mathsf{0}}=1  \label{equation}
\end{equation}%
\textit{is of the form }%
\begin{equation}
\mathcal{N}_{0,T}=c\exp \left[ \frac{TE_{\mathsf{0}}}{2}\right]
\label{value}
\end{equation}%
\textit{where }$c>0$\textit{\ is a constant depending only on }$\sigma
_{0,T} $\textit{\ and }$a_{0,T}$\textit{. Moreover, with the value (\ref%
{value}) in (\ref{fouriercoeff1}) and (\ref{fouriercoeff2}) one gets}%
\begin{equation*}
\sum_{\mathsf{n}\in \mathbb{N}^{d}}\alpha _{\mathsf{n}}\exp \left[ -TE_{%
\mathsf{n}}\right] \beta _{\mathsf{n}}=1+O\left( \exp \left[ -T\right]
\right)
\end{equation*}%
\textit{for }$T$\textit{\ sufficiently large.}

\bigskip

\textbf{Proof.} It is clear that (\ref{value}) holds because of (\ref%
{equation}) since $\hat{\alpha}_{\mathsf{0}}>0$, $\hat{\beta}_{\mathsf{0}}>0$
by virtue of the fact that the eigenfunction $\mathsf{h}_{\mathsf{0}}$
associated with the bottom of the spectrum is strictly positive in $\mathbb{R%
}^{d}$. Then, the proof that the remaining term satisfies%
\begin{equation*}
\mathcal{N}_{0,T}^{2}\sum_{\mathsf{n}\in \mathbb{N}^{d},\text{ }\mathsf{%
n\neq 0}}\hat{\alpha}_{\mathsf{n}}\exp \left[ -TE_{\mathsf{n}}\right] \hat{%
\beta}_{\mathsf{n}}=O\left( \exp \left[ -T\right] \right)
\end{equation*}%
follows from the fact that the $\hat{\alpha}_{\mathsf{n}}$'s and the $\hat{%
\beta}_{\mathsf{n}}$'s are uniformly bounded in $\mathsf{n}$, and from the
summation of the underlying geometric series which is made possible thanks
to the explicit form (\ref{energy}). \ \ $\blacksquare $

\bigskip

Then, in case of Gaussian initial-final initial data in (\ref{forwardcauchy2}%
)-(\ref{backwardcauchy2}) one gets:

\bigskip

\textbf{Proposition 4.} \textit{Assume that }$\varphi _{0}$ \textit{and }$%
\psi _{T}$\textit{\ are given by (\ref{gaussian1}) and (\ref{gaussian2}),
respectively, and let }$Z_{\tau \in \left[ 0,T\right] }$\textit{\ be the
Markovian Bernstein process associated with (\ref{forwardcauchy2})-(\ref%
{backwardcauchy2}).} \textit{Then the following statements hold: }

\textit{(a)} \textit{For all }$F_{0},F_{T}\in \mathcal{B}_{d}$\textit{\ we
have}%
\begin{eqnarray}
&&\mathbb{P}_{\mu }\left( Z_{0}\in F_{0},Z_{T}\in F_{T}\right)  \notag \\
&=&\left( 4\pi ^{2}\rho _{0}\rho _{T}\right) ^{-\frac{d}{2}}\int_{F_{0}}%
\mathsf{dx}\exp \left[ -\frac{\mathsf{1}}{2\rho _{0}}\left\vert \mathsf{x-}%
\frac{\mathsf{a}_{0}}{1+\sigma _{0}}\right\vert ^{2}\right]  \notag \\
&&\times \int_{F_{T}}\mathsf{dx}\exp \left[ -\frac{1}{2\rho _{T}}\left\vert 
\mathsf{x-}\frac{\mathsf{a}_{T}}{1+\sigma _{T}}\right\vert ^{2}\right]
+O\left( \exp \left[ -T\right] \right)  \label{jointprobability}
\end{eqnarray}%
\textit{for }$T$\textit{\ sufficiently large,} \textit{where}%
\begin{equation*}
\rho _{0,T}=\frac{\sigma _{0,T}}{1+\sigma _{0,T}}.
\end{equation*}%
\textit{In particular,}%
\begin{eqnarray}
&&\mathbb{P}_{\mu }\left( Z_{0,T}\in F_{0,T}\right) =\left( 2\pi \rho
_{0,T}\right) ^{-\frac{d}{2}}\int_{F_{0}}\mathsf{dx}\exp \left[ -\frac{1}{%
2\rho _{0,T}}\left\vert \mathsf{x-}\frac{\mathsf{a}_{0,T}}{1+\sigma _{0,T}}%
\right\vert ^{2}\right]  \notag \\
&&+O\left( \exp \left[ -T\right] \right) .  \label{probabis}
\end{eqnarray}

\textit{(b) If} $\sigma _{0}=\sigma _{T}:$ $=\sigma $ \textit{and} $\mathsf{a%
}_{0}=\mathsf{a}_{T}:=\mathsf{a}$ \textit{and if the process }$Z_{\tau \in %
\left[ 0,T\right] }$ \textit{is stationary, the preceding relations reduce to%
}%
\begin{eqnarray}
&&\mathbb{P}_{\mu }\left( Z_{t}\in F\right) =\left( 2\pi \rho \right) ^{-%
\frac{d}{2}}\int_{F}\mathsf{dx}\exp \left[ -\frac{1}{2\rho }\left\vert 
\mathsf{x-}\frac{\mathsf{a}}{1+\sigma }\right\vert ^{2}\right]  \notag \\
&&+O\left( \exp \left[ -T\right] \right)  \label{probaquarto}
\end{eqnarray}%
\textit{for }$T$\textit{\ large enough, each }$t\in \left[ 0,T\right] $%
\textit{\ and every }$F\in \mathcal{B}_{d}$,\textit{\ where }$\rho =\frac{%
\sigma }{1+\sigma }$.

\bigskip

\textbf{Proof.} From (\ref{probabilitymeasure}) and (\ref{jointproba}) one
has%
\begin{equation*}
\mathbb{P}_{\mu }\left( Z_{0}\in F_{0},Z_{T}\in F_{T}\right) =\mathcal{N}%
^{2}\int_{F_{0}}\mathsf{dx}\varphi _{0}(\mathsf{x)}\int_{F_{T}}\mathsf{dy}g(%
\mathsf{x},T,\mathsf{y})\psi _{T}(\mathsf{y)}
\end{equation*}%
where $g$ is given by (\ref{mehlerbis}), that is,%
\begin{equation}
g(\mathsf{x},T,\mathsf{y})=\hat{g}(\mathsf{x},T,\mathsf{y})+\sum_{\mathsf{n}%
\in \mathbb{N}^{d},\text{ }\mathsf{n\neq 0}}\exp \left[ -TE_{\mathsf{n}}%
\right] \mathsf{h}_{\mathsf{n}}\left( \mathsf{x}\right) \mathsf{h}_{\mathsf{n%
}}\left( \mathsf{y}\right)  \label{series}
\end{equation}%
with%
\begin{equation}
\hat{g}(\mathsf{x},T,\mathsf{y})=\pi ^{-\frac{d}{2}}\exp \left[ -\frac{1}{2}%
\left( \left\vert \mathsf{x}\right\vert ^{2}+\left\vert \mathsf{y}%
\right\vert ^{2}+dT\right) \right]  \label{approximation}
\end{equation}%
according to (\ref{energy}) and (\ref{tensorproduct}) for $\mathsf{n=0}$.
One then obtains%
\begin{eqnarray}
&&\mathcal{N}^{2}\int_{F_{0}}\mathsf{dx}\varphi _{0}(\mathsf{x)}\int_{F_{T}}%
\mathsf{dy}\hat{g}(\mathsf{x},T,\mathsf{y})\psi _{T}(\mathsf{y)}  \notag \\
&=&\mathcal{N}^{2}\pi ^{-\frac{d}{2}}\exp \left[ -\frac{dT}{2}\right] \times
\int_{F_{0}}\mathsf{dx}\exp \left[ -\frac{\left\vert \mathsf{x-a}%
_{0}\right\vert ^{2}}{2\sigma _{0}}-\frac{\left\vert \mathsf{x}\right\vert
^{2}}{2}\right]  \notag \\
&&\times \int_{F_{T}}\mathsf{dy}\exp \left[ -\frac{\left\vert \mathsf{y-a}%
_{T}\right\vert ^{2}}{2\sigma _{T}}-\frac{\left\vert \mathsf{y}\right\vert
^{2}}{2}\right] ,  \label{step}
\end{eqnarray}%
and replacing $\mathcal{N}$ by $\mathcal{N}_{0,T}$ together with the
explicit evaluation of these Gaussian integrals gives the leading term in (%
\ref{jointprobability}).

It remains to show that the contribution to (\ref{jointprobability}) coming
from the second term on the right-hand side of (\ref{series}) is
exponentially small. Writing momentarily%
\begin{equation*}
\tilde{g}(\mathsf{x},T,\mathsf{y})=\sum_{\mathsf{n}\in \mathbb{N}^{d},%
\mathsf{\ n\neq 0}}\exp \left[ -TE_{\mathsf{n}}\right] \mathsf{h}_{\mathsf{n}%
}\left( \mathsf{x}\right) \mathsf{h}_{\mathsf{n}}\left( \mathsf{y}\right)
\end{equation*}%
and estimating the absolute value of $\tilde{g}$ by using (\ref{energy}) and
(\ref{cramercharlier}), one eventually gets%
\begin{equation*}
\left\vert \tilde{g}(\mathsf{x},T,\mathsf{y})\right\vert \leq c_{d}\frac{%
\exp \left[ -\frac{\left( d+2\right) T}{2}\right] }{\left( 1-\exp \left[ -T%
\right] \right) ^{d}}
\end{equation*}%
uniformly in all $\mathsf{x,y}\in \mathbb{R}^{d}$ by summing the underlying
geometric series as before, where $c_{d}$ is a positive constant depending
only on $d$. Therefore,%
\begin{eqnarray*}
&&\mathcal{N}_{0,T}^{2}\int_{F_{0}}\mathsf{dx}\varphi _{0}(\mathsf{x)}%
\int_{F_{T}}\mathsf{dy}\left\vert \tilde{g}(\mathsf{x},T,\mathsf{y}%
)\right\vert \psi _{T}(\mathsf{y)} \\
&\leq &c_{d}\mathcal{N}_{0,T}^{2}\frac{\exp \left[ -\frac{\left( d+2\right) T%
}{2}\right] }{\left( 1-\exp \left[ -T\right] \right) ^{d}}\int_{\mathbb{R}%
^{d}}\mathsf{dx}\varphi _{0}(\mathsf{x)}\int_{\mathbb{R}^{d}}\mathsf{dy}\psi
_{T}(\mathsf{y)} \\
&=&\left( 4\pi ^{2}\sigma _{0}\sigma _{T}\right) ^{\frac{d}{2}}c_{d}\mathcal{%
N}_{0,T}^{2}\frac{\exp \left[ -\frac{\left( d+2\right) T}{2}\right] }{\left(
1-\exp \left[ -T\right] \right) ^{d}}=O\left( \exp \left[ -T\right] \right)
\end{eqnarray*}%
because of (\ref{value}), as desired. Finally (\ref{jointprobability})
implies (\ref{probabis}), and also (\ref{probaquarto}) under the hypothesis
in (b) since the function $t\mapsto \mathbb{P}_{\mu }\left( Z_{t}\in
F\right) $ given by (\ref{probability}) is then independent of $t$.\ \ $%
\blacksquare $

\bigskip

\textsc{Remark.} The first term on the right-hand side of (\ref{series})
corresponds to the minimal choice $N=1$ in the Galerkin approximation (\ref%
{galerkinapprox1}), for (\ref{energy}) and (\ref{tensorproduct}) with $%
\mathsf{n}=0$ imply that (\ref{approximation}) is%
\begin{equation*}
\hat{g}(\mathsf{x},T,\mathsf{y})=\exp \left[ -TE_{\mathsf{0}}\right] \mathsf{%
h}_{\mathsf{0}}\left( \mathsf{x}\right) \mathsf{h}_{\mathsf{0}}\left( 
\mathsf{y}\right) .
\end{equation*}%
Using once more (\ref{energy}) and (\ref{tensorproduct}) with $\mathsf{n}=0$%
, the corresponding approximation (\ref{galerkinapprox2}) for $t=T$ then
reads%
\begin{equation*}
u_{N=1}(\mathsf{x,}T)=\pi ^{-\frac{d}{4}}\mathcal{N}\hat{\alpha}_{\mathsf{0}%
}\exp \left[ -\frac{1}{2}\left( \left\vert \mathsf{x}\right\vert
^{2}+dt\right) \right] ,
\end{equation*}%
so that replacing $\mathcal{N}$ by (\ref{value}) and arguing as in the above
proofs one eventually gets%
\begin{equation}
\mathbb{P}_{\mu }\left( Z_{T}\in F_{T}\right) =\int_{F_{T}}\mathsf{dx}%
u_{N=1}(\mathsf{x,}T)\psi (\mathsf{x})+O\left( \exp \left[ -T\right] \right)
\label{probaquinto}
\end{equation}%
for $T$ sufficiently large. A similar approximation procedure applies to the
backward solution, so that in the end one obtains yet another algorithm to
compute (\ref{probabis}) since $\hat{\alpha}_{\mathsf{0}}$ and $\hat{\beta}_{%
\mathsf{0}}$ can be determined explicitly in case of Gaussian initial-final
data. It would have been difficult to evaluate (\ref{probaquinto}) directly
from (\ref{probability}) given the complicated form (\ref{forwardsolutionbis}%
). As a matter of fact, the technique used also works if the data are of the
form (\ref{hat1})-(\ref{hat2}) since $\hat{\alpha}_{\mathsf{0}}$ and $\hat{%
\beta}_{\mathsf{0}}$ are then easily determined by numerical calculations.

More generally, there is an important computational issue about (\ref%
{galerkinapprox2}) and (\ref{galerkinapprox3}), namely, that of knowing how
large one has to choose $N$ as a function of the desired degree of precision
to reconstruct $u$ and $v$. As long as error terms of the form $O\left( \exp %
\left[ -T\right] \right) $ are considered satisfactory, the above
considerations show that the choice $N=1$ is sufficient. If not, larger
values of $N$ will do.

Finally, thanks to an expansion of the form (\ref{expansion}), similar
Faedo-Galerkin approximation methods may be applied to the forward-backward
solutions of (\ref{forwardcauchy1})-(\ref{backwardcauchy1}) when the
potential function satisfies Hypothesis (H) and (\ref{convergence}), or even
more general conditions, provided that precise information be available
about the spectrum $\left( E_{\mathsf{n}}\right) _{\mathsf{n}\in \mathbb{N}%
^{d}}$ and the corresponding sequence of eigenfunctions $\left( \mathsf{f}_{%
\mathsf{n}}\right) _{\mathsf{n}\in \mathbb{N}^{d}}$. The detailed results
will be published elsewhere.

\bigskip

\textbf{Acknowledgements. }I am particularly indebted to Prof. W. Petersen
and Prof. T. Rivi\`{e}re for having made several visits in Zurich
financially possible through funds from the Forschungsinstitut f\"{u}r
Mathematik of the ETHZ, where parts of this work were carried out and whose
warm hospitality I gratefully acknowledge.

\bigskip

\bigskip

\end{document}